\documentclass[11pt]{amsart}
\usepackage{amsmath}
\usepackage{amsthm}
\usepackage{graphicx}
\usepackage{hyperref}

\newtheorem{theorem}{Theorem}

\theoremstyle{remark}

\theoremstyle{definition}

  \def\e{\epsilon}

\def\vv{{\rm v}}

\def\v{{\bf a}}
\def\aa{{\bf a}}
 \def\f_H{{\bf w}}
 \def\f{{\bf f}}
 \def\a{\alpha}

\def\Z{\mathbb{Z}}
\def\g{{\bf a}}
\def\inc{\textrm{inc}}

\begin{document}
\title[Demystifying the divisibility of Kostant partition functions]{Demystifying a divisibility property of the  Kostant partition function}
\author{Karola M\'esz\'aros}
\address{
Department of Mathematics, Massachusetts Institute of Technology, Cambridge, MA 02139
}

\begin{abstract} We study a family of  identities regarding a divisibility property of the Kostant partition function which first appeared in a paper of Baldoni and Vergne. To prove the identities,  Baldoni and Vergne used techniques of residues and called the resulting divisibility property ``mysterious."  We prove these identities entirely combinatorially and provide a natural explanation  of why the divisibility occurs. We also point out several ways to generalize the identities.

\end{abstract}

\maketitle
  
  \section{Introduction}
\label{sec:in}

The objective of this paper is to provide a natural combinatorial explanation of a  divisibility property of the Kostant partition function. The question of evaluating Kostant partition functions has been subject of much interest, without a satisfactory combinatorial answer. To mention perhaps the most famous such case: it is known that $$K_{A_n^+}(1, 2, \ldots, n, -{n+1 \choose 2})=\prod_{k=1}^{n}C_k,$$ where  $C_k=\frac{1}{k+1}{2k \choose k}$ is the Catalan number, yet there is no combinatorial proof of the above identity! 

Given such lack of understanding of the evaluation of the Kostant partition function, it seems a worthy proposition to  provide a simple explanation for its certain divisibility properties. 
We  explore divisibility properties of Kostant partition functions of types $A_n$ and $C_{n+1}$,  noting that such properties in  types $B_{n+1}$ and $D_{n+1}$ are easy consequences of the type $C_{n+1}$ case. The type $A_n$ family of identities we study  first appeared in a paper by Baldoni and Vergne \cite{BV}, where the authors prove the identities using residues, and where they call the divisibility property ``mysterious."  It is our hope that the combinatorial argument we provide   successfully demystifies the divisibility property of the Kostant partition function and provides a natural explanation why things happen the way they do.

The outline of the paper is as follows. 
In Section \ref{sec:bv} we define Kostant partition functions of type $A_n$ and prove the Baldoni-Vergne identities combinatorially. Our proof also yields an affirmative answer to a  question of Stanley \cite{S} regarding a possible bijective proof of a  special case of the Baldoni-Vergne identities.
In Section \ref{sec:c} we define Kostant partition functions of type $C_{n+1}$, relate them to flows, and show  how to modify our proof for the Baldoni-Vergne identities to obtain their analogues for type $C_{n+1}$. We also point out other possible variations of these identities in the type $C_{n+1}$ case.  

\section{The Baldoni-Vergne identities}
\label{sec:bv}

Before stating  the Baldoni-Vergne identities, we need a few definitions. Throughout this section  the graphs $G$ we consider are on the vertex set $[n+1]$ with possible multiple edges, but no loops. Denote by $m_{ij}$ the multiplicity of edge $(i, j)$, $i<j$, in $G$. To each edge $(i, j)$, $i<j$,  of $G$,  associate the positive type $A_n$ root $e_i-e_j$, where $e_i$ is the $i^{th}$ standard basis vector. Let $\{\{\a_1, \ldots, \a_N\}\}$ be the multiset of vectors corresponding to the multiset of edges of $G$ as described above. Note that $N=\sum_{1\leq i<j\leq n+1} m_{ij}$.

 The {\bf Kostant partition function}  $K_G$ evaluated at the vector $\v \in \Z^{n+1}$ is defined as

\begin{equation} \label{kost} K_G(\v)= \# \{ (b_{i})_{i \in [N]} \mid \sum_{i \in [N]} b_{i}  \a_i =\v \textrm{ and } b_{i} \in \Z_{\geq 0}\}.\end{equation}

That is, $K_G(\v)$ is the number of ways to write the vector $\v$ as a nonnegative linear combination of the positive type $A_n$ roots corresponding to the edges of $G$, without regard to order.  Note that in order for $K_G(\v)$ to be nonzero, the partial sums of the coordinates of $\v$ have to satisfy $v_1+\ldots+v_i \geq 0$, $i\in [n]$, and $v_1+\ldots+v_{n+1}=0$.

We now proceed to state and prove Theorem \ref{thm:bv} which first appeared in  \cite{BV}.   Baldoni and Vergne  gave a proof of it using residues, and called the result ``mysterious."  We provide a natural combinatorial explanation of the result. Our explanation also answers a question of Stanley in affirmative, which he posed in \cite{S}, regarding a possible bijective proof of a  special case of the Baldoni-Vergne identities.

For brevity, we write $G-e$, or $G-\{e_1, \ldots, e_k\}$, to mean a graph obtained from  $G$ with edge $e$, or edges $e_1, \ldots, e_k$, deleted. 

\begin{theorem} \label{thm:bv} (\cite{BV})
Given a connected graph $G$ on the vertex set $[n+1]$ with $m_{n-1, n}=m_{n-1, n+1}=m_{n, n+1}=1$ and such that $$\frac{m_{j, n-1}+m_{j, n}+m_{j, n+1}}{m_{j, n-1}}=c,$$ 

\noindent for some constant $c$ independent of $j$ for $j \in [n-2]$, we have that for any $\aa=(a_1, \ldots, a_n, -\sum_{i=1}^n a_i) \in \Z^{n+1}$ 

\begin{equation} \label{c} K_G(\aa)=(\frac{a_1+ \cdots+ a_{n-2}}{c}+ a_{n-1}+1)   K_{G- (n-1, n)}(\aa).\end{equation}
\end{theorem}

\bigskip

Before proceeding to the formal proof of Theorem \ref{thm:bv} we outline it, to fully expose the underlying combinatorics. Rephrasing equation \eqref{kost}, $K_G(\v)$  counts the number of {\bf flows} $\f_G=(b_i)_{i\in N}$ on graph $G$ satisfying $$\sum_{i \in [N]} b_{i}   \a_i =\v \textrm{ and } b_{i} \in \Z_{\geq 0}.$$ In the proof of Theorem \ref{thm:bv} we introduce the concept of {\bf partial flows} $\f_H$   and the following are the key statements we prove:
\begin{itemize}
\item The elements of the set of  partial flows are in bijection with the flows on $G-(n-1, n)$ that the Kostant partition function $K_{G-(n-1, n)}(\aa)$ counts. That is, $$\#\text{ partial flows}=K_{G-(n-1, n)}(\aa).$$

\item The elements of the  multiset of partial flows  $\f_H$, where the cardinality of the multiset is  $\frac{a_1+ \cdots+ a_{n-2}}{c}+ a_{n-1}+1$ times the cardinality of the set of partial flows, are  in bijection with the flows on $G$ that the Kostant partition function $K_{G}(\aa)$ counts. That is, $$(\frac{a_1+ \cdots+ a_{n-2}}{c}+ a_{n-1}+1)  \#\text{ partial flows}=K_{G}(\aa).$$

\end{itemize}

The above two statements imply a bijection between the elements of the multiset of  flows on $G-(n-1, n)$ that the Kostant partition function $K_{G-(n-1, n)}(\aa)$ counts,  where the cardinality of the multiset is  $\frac{a_1+ \cdots+ a_{n-2}}{c}+ a_{n-1}+1$ times the cardinality of the set of  flows counted by  $K_{G-(n-1, n)}(\aa)$, and  the flows on $G$ that the Kostant partition function $K_{G}(\aa)$  counts, yielding 

$$K_G(\aa)=(\frac{a_1+ \cdots+ a_{n-2}}{c}+ a_{n-1}+1)   K_{G- (n-1, n)}(\aa).$$
\bigskip

We now proceed to the formal proof of Theorem \ref{thm:bv}. 

\bigskip

\noindent \textit{Proof of Theorem \ref{thm:bv}}.
Let $\{\{\a_1, \ldots, \a_N\}\}$ be the multiset of vectors corresponding to the edges of $G$. Let $\a_N=e_{n-1}-e_n, \a_{N-1}=e_{n-1}-e_{n+1},$ and $\a_{N-2}=e_{n}-e_{n+1}.$ Then equation \eqref{c} can be rewritten as

\begin{align} \label{c1} & \# \{ (b_i)_{i \in [N]} \mid \sum_{i=1}^N b_i   \a_i =\aa\} = \nonumber \\
&(\frac{a_1+ \cdots+ a_{n-2}}{c}+ a_{n-1}+1)     \# \{ (b_i)_{i \in [N-1]} \mid \sum_{i=1}^{N-1} b_i   \a_i =\aa\}.
\end{align}

Consider a flow $\f_H=(b_i)_{i \in [N-3]}$, $b_i \in \Z_{\geq 0}$,  of the edges of the graph $H:=G-\{(n-1, n), (n-1, n+1), (n, n+1)\}$. We call $\f_H$ {\bf partial} if  $$\sum_{i=1}^{N-3} b_i   \a_i =(a_1, \ldots, a_{n-2}, x_1, x_2, x_3),$$ \noindent for some $x_1, x_2, x_3 \in \Z$. 

Notice that given a partial flow $\f_H=(b_i)_{i \in [N-3]}$, it can be extended uniquely to a flow $\f_{G-\{(n-1, n)\}}=(b_i)_{i \in [N-1]}$,  $b_i \in \Z_{\geq 0}$, on $G-\{(n-1, n)\}$ such that 
$\sum_{i=1}^{N-1} b_i   \a_i =\aa$. Furthermore, this correspondence is a bijection. Therefore, 

\begin{equation} \label{H} \# \{ (b_i)_{i \in [N-1]} \mid \sum_{i=1}^{N-1} b_i   \a_i =\aa\}=\sum_{\f_H} 1,\end{equation} 

\noindent where the summation runs over all partial flows $\f_H$.

Also,  given a partial flow $\f_H$ with  $Y_{i}(\f_H)$, $i \in \{n-1, n, n+1\}$, denoting the total {\bf inflow} into vertex $i \in \{n-1, n, n+1\}$ in $H$, that is the sum of all the flows $b_i$ on edges of $H$ incident to $i  \in \{n-1, n, n+1\}$,  the partial flow $\f_H$ can be extended in $Y_{n-1}(\f_H)+a_{n-1}+1$ ways to a flow $\f_{G}=(b_i)_{i \in [N]}$,  $b_i \in \Z_{\geq 0}$, of $G$ such that 
$\sum_{i=1}^{N} b_i   \a_i =\aa$. Furthermore, given a flow $\f_G=(b_i)_{i \in [N]}$,  $b_i \in \Z_{\geq 0}$ such that  $\sum_{i=1}^{N} b_i   \a_i =\aa$, there is a unique partial flow $\f_H=(b_i)_{i \in [N-3]}$ from which it can be obtained.  Therefore, 

\begin{equation} \label{G}  \# \{ (b_i)_{i \in [N]} \mid \sum_{i=1}^N b_i   \a_i =\aa \} =\sum_{\f_H} (Y_{n-1}(\f_H)+a_{n-1}+1),\end{equation} 

\noindent where the summation runs over all partial flows $\f_H$.

Note that since $$\frac{m_{j, n-1}+m_{j, r}+m_{j, n+1}}{m_{j, n-1}}=c,$$ 

\noindent for some constant $c$ independent of $j$ for $j \in [n-2]$, it follows that 

$$c   \sum_{\f_H} Y_{n-1}(\f_H)=\sum_{\f_H} (Y_{n-1}(\f_H)+Y_{n}(\f_H)+Y_{ n+1}(\f_H))=\sum_{\f_H} (a_1+\cdots+a_{n-2}),$$

that is 

\begin{equation}\label{eq:y} \sum_{\f_H} Y_{n-1}(\f_H)=\sum_{\f_H} \frac{(a_1+\cdots+a_{n-2})}{c}.\end{equation}

Thus, equation \eqref{G} can be rewritten as

\begin{align} \label{G1} &  \# \{ (b_i)_{i \in [N]} \mid \sum_{i=1}^N b_i   \a_i =\aa \} =\sum_{\f_H} (\frac{(a_1+\cdots+a_{n-2})}{c}+a_{n-1}+1) \\ &=(\frac{(a_1+\cdots+a_{n-2})}{c}+a_{n-1}+1)    \sum_{\f_H}  1\\ &=(\frac{(a_1+\cdots+a_{n-2})}{c}+a_{n-1}+1)   \# \{ (b_i)_{i \in [N-1]} \mid \sum_{i=1}^{N-1} b_i   \a_i =\aa\}, \end{align} 

\noindent where the first equality uses equations \eqref{G} and \eqref{eq:y}, and the third equality uses equation \eqref{H}. 
\qed

 \section{Type $C_{n+1}$ Kostant partition functions and the Baldoni-Vergne identities}
\label{sec:c}

We now show two generalizations of Theorem \ref{thm:bv} in the type $C_{n+1}$ case. We first  give the necessary definitions and explain the notion of flow in the context of signed graphs. 
Throughout this section  the graphs $G$ on the vertex set $[n+1]$ we consider are signed, that is there is a sign $\e \in \{+, -\}$ assigned to each of its edges, with possible multiple edges, and all loops labeled positive. Denote by $(i, j, -)$ and $(i, j, +)$, $i \leq j$, a negative and a positive edge, respectively. Denote by $m_{ij}^\e$ the multiplicity of edge $(i, j, \e)$ in $G$, $i\leq j$, $\e \in \{+, -\}$. 
To each edge $(i, j, \e)$, $i\leq j$,  of $G$,  associate the positive type $C_{n+1}$ root $\vv(i,j, \e)$, where $\vv(i,j, -)=e_i-e_j$ and $\vv(i,j, +)=e_i+e_j$. Let $\{\{\a_1, \ldots, \a_N\}\}$ be the multiset of vectors corresponding to the multiset of edges of $G$ as described above. Note that $N=\sum_{1\leq i<j\leq n+1} (m_{ij}^-+m_{ij}^+)$.

 The {\bf Kostant partition function}  $K_G$ evaluated at the vector $\v \in \Z^{n+1}$ is defined as

$$K_G(\v)= \# \{ (b_{i})_{i \in [N]} \mid \sum_{i \in [N]} b_{i}   \a_i =\v \textrm{ and } b_{i} \in \Z_{\geq 0}\}.$$

That is, $K_G(\v)$ is the number of ways to write the vector $\v$ as a nonnegative linear combination of the positive type $C_{n+1}$ roots corresponding to the edges of $G$, without regard to order.  

Just like in the type $A_n$ case, we would like to think of the vector $(b_{i})_{i \in [N]} $ as a {\bf flow}. For this we here give a precise definition of flows in the type $C_{n+1}$ case, of which type $A_{n}$ is of course a special case.   

Let $G$ be a signed graph on the vertex set $[n+1]$. Let $\{\{e_1, \ldots, e_N\}\}$ be the multiset of edges of $G$, and $\{\{\a_1, \ldots, \a_N\}\}$  the multiset of vectors corresponding to the multiset of edges of $G$. Fix an integer  vector $\g=(a_1, \ldots, a_n, a_{n+1}) \in \Z^{n+1}$.  A {\bf nonnegative integer} \textbf{$\g$-flow} $\f_G$ on $G$ is a vector $\f_G=(b_i)_{i \in [N]}$, $b_i \in \Z_{\geq 0}$  such that for all $1\leq i \leq n+1$, we have 

 \begin{equation} \label{eqn:flow}
 \sum_{e \in E, \inc(e, v )=-} b(e)+a_v= \sum_{e \in E, \inc(e, v)=+} b(e)+\sum_{e=(v, v, + )} b(e),
  \end{equation}
  
\noindent   where $b(e_i)=b_i$, $\inc(e, v)=-$ if edge $e=(g, v, -)$, $g <v$, and $\inc(e, v)=+$ if $e=(g, v, +)$, $g <v$, or $e=( v, j, \e)$, $v <j,$ and $\e \in \{+, -\}$.

  Call $b(e)$  the \textbf{flow} assigned to edge $e$ of $G$. If the edge $e$ is negative, one can think of $b(e)$ units of fluid flowing on $e$ from its smaller to its bigger vertex. If the edge $e$ is positive, then one can think of $b(e)$ units of fluid flowing away both from $e$'s smaller and bigger vertex to infinity. Edge $e$ is then a ``leak" taking away $2b(e)$ units of fluid.

From the above explanation it is clear that if we are given an $\aa$-flow $\f_G$ such that  $\sum_{e=(i, j, +), i \leq j}b(e)=y$, then $\aa=(a_1, \ldots, a_n, 2y-\sum_{i=1}^n a_i)$.   It is then a matter of checking the defintions to see that for a signed graph $G$ on the vertex set $[n+1]$ and vector  $\g=(a_1, \ldots, a_n, 2y-\sum_{i=1}^n a_i) \in \Z^{n+1}$,  the number of nonnegative integer $\g$-flows on $G$ is equal to $K_G( \g)$. 

Thinking of $K_G(\g)$ as the number of nonnegative integer $\g$-flows on $G$, there is a straightforward generalization of Theorem \ref{thm:bv} in the type $C_{n+1}$ case: 

\begin{theorem} \label{thm:bv2} 
Given a connected signed graph $G$ on the vertex set $[n+1]$ with $m_{n-1, n}^-=m_{n-1, n+1}^-=m_{n, n+1}^-=1$, $m_{j, n-1}^+=m_{j, n}^+=m_{j, n+1}^+=0$, for $j \in [n+1]$,  and such that $$\frac{m_{j, n-1}^-+m_{j, n}^-+m_{j, n+1}^-}{m_{j, n-1}^-}=c,$$ 

\noindent for some constant $c$ independent of $j$ for $j \in [n-2]$, we have that for any $\aa=(a_1, \ldots, a_n, 2y-\sum_{i=1}^n a_i) \in \Z^{n+1}$, 

\begin{equation} \label{cC} K_G(\aa)=(\frac{a_1+ \cdots+ a_{n-2}-2y}{c}+ a_{n-1}+1)   K_{G- (n-1, n)}(\aa).\end{equation}
\end{theorem}

\bigskip

The proof of Theorem \ref{thm:bv2} proceeds analogously to that of Theorem \ref{thm:bv}. Namely,  define \textbf{partial flows} $\f_H=(b_i)_{i\in [N-3]}$ on $H:=G-\{(n-1, n, -), (n-1, n+1, -), (n, n+1, -)\}$ such that  $$\sum_{i=1}^{N-3} b_i   \a_i =(a_1, \ldots, a_{n-2}, x_1, x_2, x_3),$$ \noindent for some $x_1, x_2, x_3 \in \Z$ and the sum of flows on positive edges is $y$. 

Then, one can prove:
\begin{itemize}
\item The elements of the set partial flows are in bijection with the nonnegative integer $\aa$-flows on $G-(n-1, n)$. That is, $$\#\text{ partial flows}=K_{G-(n-1, n)}(\aa).$$

\item The elements of the  multiset of partial flows  $\f_H$, where the cardinality of the multiset is  $\frac{a_1+ \cdots+ a_{n-2}-2y}{c}+ a_{n-1}+1$ times the cardinality of the set of partial flows, are  in bijection with the nonnegative integer $\aa$-flows on $G$ . That is, $$(\frac{a_1+ \cdots+ a_{n-2}-2y}{c}+ a_{n-1}+1)  \#\text{ partial flows}=K_{G}(\aa).$$

\end{itemize}

The above two statements imply a bijection between the elements of the multiset of  nonnegative integer $\aa$-flows on $G-(n-1, n)$, where the cardinality of the multiset is   $\frac{a_1+ \cdots+ a_{n-2}-2y}{c}+ a_{n-1}+1$ times the cardinality of the set of nonnegative integer $\aa$-flows on $G-(n-1, n)$, and  the nonnegative integer $\aa$-flows on $G$, yielding 

$$K_G(\aa)=(\frac{a_1+ \cdots+ a_{n-2}-2y}{c}+ a_{n-1}+1)   K_{G- (n-1, n)}(\aa).$$
\bigskip

Note that the requirement that only negative edges are incident to the vertices $n-1, n, n+1$ in $G$ stems from the fact that we need to make sure, in order for our counting  arguments from the proof of Theorem \ref{thm:bv} to work, that we can always assign nonnegative flows to the edges $(n-1, n, -), (n-1, n+1, -), (n, n+1, -)$ and also, that in case we are extending a partial flow $\f_H$ to a flow on $G$, we can extend it in $Y_{n-1}(\f_H)+a_{n-1}+1$ ways. These properties will be satisfied, if we insure that ``inflows" at the vertices 
$n-1$ and $n$ are at least $-a_{n-1}$ and $-a_n$, respectively. To simplify the formulation, we will assume that there are no loops at the vertices $n-1, n, n+1$, though the following theorem could also be adopted to a somewhat more general setting.

\begin{theorem} \label{thm:bv3} 
Given a connected signed graph $G$ on the vertex set $[n+1]$ with $m_{n-1, n}^-=m_{n-1, n+1}^-=m_{n, n+1}^-=1$, $m_{i,j}^+=0$, for $i, j \in \{n-1, n, n+1\}$,  and such that $$\frac{m_{j, n-1}^\e+m_{j, r}^\e+m_{j, n+1}^\e}{m_{j, n-1}^\e}=c,$$ 

\noindent for $\e \in \{+, -\}$ and for some constant $c$ independent of $j$ for $j \in [n-2]$, we have that for any $\aa=(a_1, \ldots, a_n, 2y-\sum_{i=1}^n a_i) \in \Z^{n+1}$,  $y \leq a_{n-1}+1, a_n+1$,

\begin{equation} \label{cC} K_G(\aa)=(\frac{a_1+ \cdots+ a_{n-2}-2y}{c}+ a_{n-1}+1)   K_{G- (n-1, n)}(\aa).\end{equation}
\end{theorem}

The proof technique of Theorem \ref{thm:bv3} is analogous to that of Theorem \ref{thm:bv}. We invite the reader to check each step of the proof of Theorem \ref{thm:bv} and see how they can be adapted to prove Theorem \ref{thm:bv3}.

\section*{Acknowledgement} I  thank Richard Stanley for his intriguing slides and for pointing out the work of Baldoni and Vergne. I also thank   Alejandro Morales for numerous conversations about flows and the Kostant partition function.


\begin{thebibliography}{99}



\bibitem[BV]{BV} {\sc W. ~Baldoni, M. ~Vergne},
Kostant partition functions and flow polytopes, 
in {\it  Transformation Groups 13,} Birkh\"auser, Boston,   2008, 447-469.




\bibitem[S]{S} {\sc R. ~Stanley},
 {\it Acyclic flow polytopes and Kostant's partition function,} Conference transparencies 2000,  {\tt http://math.mit.edu/$\sim$rstan/trans.html}.

\end{thebibliography}
\end{document}